%% file: Preprint_2_DoldKan.tex
\documentclass[12pt]{article}

\usepackage{erics_preprints}

\input xy
\xyoption{all}
\xyoption{poly}
\xyoption{matrix}

\begin{document}

\title{Nonabelian Dold-Kan Decompositions for Simplicial and Symmetric-Simplicial Groups}
\author{Eric Ram\'on Antokoletz}
\date{}
\maketitle

\begin{abstract}

We extend the nonabelian Dold-Kan decomposition for simplicial groups of
Carrasco and Cegarra \cite{CarrascoCegarra} in two ways. First, we show that the
total order of the subgroups in their decomposition belongs to a family
of total orders all giving rise to Dold-Kan decompositions. We exhibit a
particular partial order such that the family is characterized as
consisting of all total orders extending the partial order. Second, we
consider symmetric-simplicial groups and show that, by using a specially
chosen presentation of the category of symmetric-simplicial operators,
new Dold-Kan decompositions exist which are algebraically much simpler
than those of \cite{CarrascoCegarra} in the sense that the commutator of
two component subgroups lies in a single component subgroup.

\end{abstract}

\tableofcontents
		
\setcounter{secnumdepth}{2}
\setcounter{section}{0}

\pagebreak
							
	\section{Introduction}						\input{Intro_2_DoldKan}

	\section{Preliminaries}						
		\subsection{The Simplicial and Symmetric-Simplicial Categories}			\input{SymmMooreComp}

		\subsection{Moore Complexes of Simplicial and Symmetric-Simplicial Groups}	\input{Moorecomplex}
	\section{Dold-Kan Decompositions for Simplicial Groups}					\input{DoldKanDecomps}

	\section{Simplicial Identities and the Length-Product Partial Order}			\input{ProofLemmaIdsNew}

	\section{On the Binary Total Order}							\input{BinaryTO}
	\section{Alternative Dold-Kan Decompositions for Symmetric Simplicial Groups}		\input{SymmDoldKan}

\bibliographystyle{amsalpha}
\bibliography{Preprint_2_DoldKan}

\end{document}

%% file: Intro_2_DoldKan.tex
	\label{Intro_SDPs}


In this introduction we refer to the simplicial category $\ord$
and the notion of simplicial object in a category $\C$, as well as the
Moore complex $M(G)$ of a simplicial group $G$. These
notions are reviewed in the two
subsections of section \ref{Moorecomplex}.

The classical Dold-Kan theorem (\cite{May}, \cite{GoerssJardine}) says that,
in the case in which $G$ is a simplicial abelian group $A$,
the terms $A_n$ of $A$ decompose as direct sums of copies terms $M_n = M_n(A)$ of its
Moore complex as exemplified below for the cases $n=0,1,2,3$.
	\begin{align*}
		A_0 &\cong M_0 \\
		A_1 &\cong M_1 \oplus s_0M_0 \\
		A_2 &\cong M_2 \oplus s_0M_1 \oplus s_1M_1 \oplus s_1s_0M_0 \\
		A_3 &\cong M_3 \oplus s_0M_2 \oplus s_1M_2 \oplus s_2M_2 \oplus s_1s_0M_1 
				 \oplus s_2s_0M_1 \oplus s_2s_1M_1 \oplus s_2s_1s_0M_0
	\end{align*}
For general $n\ge0$ one has the isomorphism (\cite{GoerssJardine})
	\begin{align*}
		A_n &\cong M_n\oplus\left(\bigoplus_{k=0}^{n-1}\bigoplus_{\substack{\sigma\in\ord \\ \sigma:[n]\tra[k]}}{\sigma^*M_{k}}\right)
	\end{align*}
which says that the summands of $A_n$ are precisely $M_n$ together with
the images in $A_n$ of lower terms of the Moore complex under all
possible degeneracy operators $\sigma^*$ landing in $A_n$. Using the
formal structure of this formula, the terms of $A$ and the action of the
simplicial operators on $A$ can be completely reconstructed from the
data contained in $M(A)$. 

In this way, the Dold-Kan theorem shows that the functor $M$ is an
equivalence of categories between the category ${\bf SAb}$ of simplicial
abelian groups and the category ${\bf Ch}_+({\mathbb Z})$ of nonnegative
chain complexes of abelian groups. Under this equivalence, referred to
as the {\it Dold-Kan correspondence}, homomorphic homotopy equivalences
of simplicial abelian groups correspond to quasi-isomorphisms of chain
complexes. In this sense the Dold-Kan correspondence may be said, on the one
hand, to elevate the homological algebra of ${\bf Ch}_+({\mathbb Z})$ to
the level of homotopy theory, and on the other hand, to reduce the
homotopy theory of ${\bf SAb}$ to homological algebra.

From this perspective, the general nonabelian case is especially
interesting because of the classical result that the category of
simplicial groups possesses a homotopy theory equivalent to that of
pointed, connected topological spaces (see \cite{KanLine} and
\cite{Quillen}). It would then seem that a nonabelian Dold-Kan theorem
could provide a way of translating the homotopy theory of spaces---a
subject known for the difficulty in calculating its basic
invariants---to a context resembling homological algebra, in which many
computational techniques have been developed and standardized. Indeed
such a theorem was provided by P. Carrasco and A. M. Cegarra, and they
demonstrated its bridging role by using it to describe algebraic models
for homotopy types called {\it hypercrossed complexes} and giving
explicit descriptions of truncated hypercrossed complexes yielding
models for homotopy 3-types (see \cite{CarrascoCegarra}).

In order to do this, Carrasco and Cegarra first identified the
appropriate notion of {\it higher semidirect product (SDP)} which takes
the place of the direct sums in the classical Dold-Kan decomposition
(see Definition \ref{iSDP} below or the paper \cite{preprint1} for a
more in-depth investigation of SDPs). With this notion in hand, they
prove that $G_n$ is an SDP of copies of terms $M_n := M_n(G)$ of its
Moore complex, as shown here for the cases $n=0,1,2,3$.
	\begin{align*}
		G_0 &= M_0 \\
		G_1 &= M_1 \rtimes s_0M_0 \\
		G_2 &= M_2 \rtimes s_0M_1 \rtimes s_1M_1 \rtimes s_1s_0M_0 \\
		G_3 &= M_3 \rtimes s_0M_2 \rtimes s_1M_2 \rtimes s_1s_0M_1 \rtimes s_2M_2 \rtimes s_2s_0M_1 \rtimes s_2s_1M_1 \rtimes s_2s_1s_0M_0
	\end{align*}
One feature of SDPs is that the order of the factors in the
decomposition is an essential part of the structure. Therefore an
important aspect of the result of \cite{CarrascoCegarra} is a
total ordering of the subgroups $s_{i_k}\ldots s_{i_1} M_{n-k} \subseteq
G_n$ giving rise to the SDP decomposition of $G_n$.

The first goal of the present paper, accomplished in sections
\ref{DoldKanDecomps} and \ref{ProofLemmaIdsNew}, is to extend Carrasco
and Cegarra's nonabelian Dold-Kan decomposition in the following way.
Under a particular choice of convention for the Moore complex, their
total order takes the form exemplified above, so that we are justified
in referring to it as the {\it binary order}. We exhibit a special
partial order on the same collection of subgroups, the {\it
length-product partial order}, such that any total order respecting the
length-product partial order will also yield an SDP decomposition of
$G_n$.

The data necessary to describe hypercrossed complexes, consisting
essentially of components of commutator brackets, is admittedly rather
complicated. Since it follows from results of Dwyer-Hopkins-Kan (see
\cite{DwyerHopkinsKan}) and the author (to appear elsewhere) that the
category of {\it symmetric-simplicial groups} (for a definition see
paragraph just after Remark \ref{nonredundant}) also possesses a
homotopy theory equivalent to that of pointed connected topological
spaces, one might wonder whether their extra structure might give rise
to hypercrossed complexes of a more manageable character. 

We take up this question in section \ref{SymmDoldKan}, relying on a
particular presentation of the category of symmetric-simplicial
operators (see Theorem \ref{altpresoffin}) derived in \cite{preprint4}
especially for this purpose. The answer is that, if $G$ is a
symmetric-simplicial group, there are new Dold-Kan-type decompositions
available for it in addition to the ones described in section
\ref{DoldKanDecomps}. For these new decompositions, many more orderings
of the components are available than in the earlier case. In effect, the
length-product partial order is replaced by the partial order given by
inclusion (of sets of indices), and requiring a total order to extend
this partial order places many fewer constraints on it. This extra
flexibility is a reflection of the fact that the commutators coming from
pairs of the component subgroups thus obtained have only a single
nontrivial component, giving the data constituting {\it symmetric
hypercrossed complexes} a sleeker form than that of the original
hypercrossed complexes. This will be demonstrated in a forthcoming
publication.

%% file: SymmMooreComp.tex

\label{SymmMooreComp}

Recall that a simplicial object in a category $\C$ is defined as a
functor $\ord^{op}\ra\C$ where $\ord$ is the category of finite ordered
sets $[n] := \{0,1,\ldots,n\}$ and (not necessarily strictly) monotonic
maps between them. $\ord$ is generated by the {\it cofaces} $d_i$ and
{\it codegeneracies} $s_i$ which are defined by saying that
$d_i:[n]\ra[n+1]$ is the unique monotonic injection for which every
fiber (i.e., preimage of a singleton) has one element except for the
fiber of $i$ which has none, while $s_i:[n]\ra[n-1]$ is the unique
monotonic surjection for which every fiber has one element except for
the fiber of $i$ which has two. Here we recall the well-known
presentation of the category $\ord$ via generators and relations (see
\cite{MacLane} for an elegant proof).  Since in this paper
we will follow the tradition of applying simplicial operators on the {\it left}
of simplicial objects, we state the relations in opposite form, i.e., as
a presentation of $\ord^{op}$.

\pagebreak

	$$\text{{\bf \underline{The Simplicial Identities.}}}$$
	$$\begin{array}{ll}
		{d_id_j = \left\{\begin{array}{cl}
				d_{j-1}d_{i} & ~~\textrm{if}~~ i < j \\
				d_jd_{i+1} & ~~\textrm{if}~~ i \ge j  
			\end{array}\right.} &
		  \\[12pt]
		{s_is_j = \left\{\begin{array}{cl}
				s_{j+1}s_{i} & ~~\textrm{if}~~ i \le j \\
				s_js_{i-1} & ~~\textrm{if}~~  i > j 
			\end{array}\right.} & \\[13pt]
		{d_is_j = \left\{\begin{array}{cl}
				s_{j-1}d_{i} & ~~\textrm{if}~~ i < j \\
				\text{id} & ~~\textrm{if}~~ i = j \mathrm{~or~} j+1 \\
				s_jd_{i-1} & ~~\textrm{if}~~  i \ge j+2
			\end{array}\right.} &
			s_id_j = \left\{\begin{array}{cl}
				d_{j+1}s_i & ~~\textrm{if}~~ i < j \\
				d_js_{i+1} & ~~\textrm{if}~~ i \ge j 
			\end{array}\right. \\
\end{array}$$

\rmk{\label{nonredundant}
	These identities are usually written in a nonredundant form. Here,
	and in all other presentations below, we have included all possible
	situations that arise when interchanging two generators, thus incurring
	a certain amount of redundancy. In particular, the table above (as well
	as the others to come) is arranged so that all identities in the right
	column follow from the identities to their left, so that all identities
	in the right column are redundant. Nevertheless, some redundancies
	remain in the left column. 	
}

Similarly as above, a {\it symmetric-simplicial object} in a category $\C$ is
defined as a functor $\fin^{op}\ra\C$ where $\fin$ is the category of
finite ordered sets $[n] := \{0,1,\ldots,n\}$ and {\it all} maps between
them. $\fin$ is generated by its subcategory $\ord$ together with the
groups $\sym_n$ for all $n\ge0$ of all permutations of the set $[n]$ (so
$\sym_n$ is a symmetric group on $n+1$ elements). A presentation of
$\fin$ has been given by Marco Grandis (see \cite{Grandis}) using the
generators of $\ord$ together with the transpositions $t_i\in\sym_n$ (see
Definition \ref{adjtransdef} below). Here we shall need an alternative
presentation of $\fin$ derived in \cite{preprint4}, which uses the
generators $d_i$ and $t_i$ but replaces the codegeneracies $s_i$ with
the {\it quasi-codegeneracies} $u_i$ (see Definition \ref{quasicodegens}
below).  In order to state the alternative presentation, we first define
the relevant maps explicitly.

\defn{ \label{adjtransdef}
	The following maps in $\fin$ are called the 
	{\it adjacent transpositions} and are defined as follows.
	$$t_i = t_i^{(n)}: [n] \lra [n] ~ \mfor n \ge 1 \mand 0 \le i \le n-1$$
	$$k \mapsto \left\{\begin{array} {cl}
				k & \mfor k \neq i,~i+1 \\
				i+1 & \mfor k = i \\
				i & \mfor k = i+1
			\end{array}\right.$$
}

\defn{
	The following maps in $\fin$ are called the {\it standard cyclic
	permutations}. 
		$$z_i = z_i^{(n)}: [n] \lra [n] ~ \mfor n \ge 1 \mand 0 \le i \le n$$ 
		$$k \mapsto \left\{\begin{array} {cl}
				k+1 & \mfor 0 \le k \le i-1 \\
				0 & \mfor k = i \\
				k & \mfor k > i
			\end{array}\right.$$
	Note that $z_i$ is an $(i+1)$-cycle on the elements $0,1,\ldots, i$. In
	particular, $z_0$ is the identity. One may equivalently take the
	following formula in $\fin^{op}$ as a definition of the corresponding
	symmetric-simplicial operator $z_i$ for $n \ge 0 \mand 0 \le i \le n$.
		$$z_i = z_i^{(n)} := t_{i-1} \ldots t_1 t_0$$ 
}

\defn{ \label{quasicodegens}
	The following maps in $\fin$ will be referred to as the {\it elementary
	quasi-codegeneracies}.
	$$u_i = u_i^{(n)}: [n+1] \lra [n] ~ \mfor n \ge 0 \mand 1 \le i \le n+1$$
	$$k \mapsto \left\{\begin{array} {cl}
				0 & \mfor k = 0 \textrm{~or~} i \\
				k & \mfor 1 \le k \le i-1 \\
				k-1 & \mfor k > i
			\end{array}\right.$$
	In particular, $u_1$ coincides with $s_0$. Note $u_0$ is not defined.
	One may equivalently define the {\it elementary quasi-degeneracy
	operators} $u_i\in\fin^{op}$ in terms of the $s_i$ and $z_i$ by means of
	the following formula holding in $\fin^{op}$ for $i \ge 1$. 
		$$u_i := z_{i-1}^{-1}s_{i-1}z_{i-1}$$
}

\smallskip

Just as for simplicial objects above, we will also apply
symmetric-simplicial operators on the {\it left} of symmetric-simplicial
objects, so we state the alternative presentation of $\fin$ in opposite
form, i.e., as a presentation of $\fin^{op}$.

\pagebreak
\begin{thm}\label{QuasiIds}\label{altpresoffin}
	The generators $d_i$, $u_i$, and $t_i$ together with the 
	following relations constitute a presentation of $\fin^{op}$.
	$$\begin{array}{ll}
		{d_id_j = \left\{\begin{array}{cl}
				d_{j-1}d_{i} & ~~\mif~~ i < j \\
				d_jd_{i+1} & ~~\mif~~ i \ge j  
			\end{array}\right.} & \\[15pt]
		{d_iu_j = \left\{\begin{array}{cl}
				z_{j-1} & \phantom{0 \neq}~~~\mif~~ i=0 \\
				u_{j-1}d_{i} & ~~\mif~~          0 \neq  i < j \\
				\id & \phantom{0 \neq}~~~\mif~~  i = j \\
				u_jd_{i-1}   & \phantom{0 \neq}~~~\mif~~  i > j
			\end{array}\right.} &
		{u_id_j = \left\{\begin{array}{ll}
				d_{j+1}u_i & ~~\mif~~ i \le j \\
				d_ju_{i+1} & ~~\mif~~ i \ge j \neq 0 \\
				d_1u_{i+1}t_0 & ~~\mif~~ j = 0
			\end{array}\right.} \\[27pt]
		{u_iu_j = \left\{\begin{array}{cl}
				u_{j+1}u_{i} & ~~\mif~~ i \le j \\
				u_ju_{i-1} & ~~\mif~~ i > j
			\end{array}\right.} \\[13pt]
		{t_it_j = \left\{\begin{array}{cl}
			\id & ~~\mif~~ i = j \\
			t_jt_i & ~~\mif~~ \abs{i-j} \ge 2 \\
			(t_jt_i)^2 & ~~\mif~~ \abs{i-j} = 1
			\end{array}\right.} \\[20pt]
		{d_it_j = \left\{\begin{array}{cl}
			t_{j-1}d_{i} & ~~\mif~~ i < j \\
			d_{i+1} & ~~\mif~~ i = j \\
			d_{i-1} & ~~\mif~~ i = j+1 \\
			t_jd_i & ~~\mif~~ i \ge j+2
			\end{array}\right.}  &
		{t_id_j = \left\{\begin{array}{cl}
			d_jt_i & ~~\mif~~ i \le j-2 \\
			d_jt_{i+1}t_it_{i+1} & ~~\mif~~ i = j-1 \\
			d_jt_{i+1} & ~~\mif~~ i \ge j
			\end{array}\right.} \\[20pt]
		{t_iu_j = \left\{\begin{array}{cl}
				u_jt_i & ~~\mif~~ 0 \neq i \le j-2 \\
				u_{j-1} & ~~\mif~~ 0 \neq i = j-1 \\
				u_{j+1} & \phantom{0 \neq}~~~\mif~~ i = j \\
				u_jt_{i-1} & \phantom{0 \neq}~~~\mif~~ i > j \\
			\end{array}\right.} &
		{u_it_j = \left\{\begin{array}{cl}
				t_{j+1}u_i & ~~\mif~~ i \le j \\
				t_jt_{j+1}u_{i-1} & \begin{array}{r} ~\mif~~ i = j+1 \\ \mand j \neq 0\end{array} \\
				t_ju_i & \begin{array}{r} ~\mif~~ i \ge j+2 \\ \mand j \neq 0\end{array}
			\end{array}\right.} \\
		{t_0u_1 = u_1} \\
		{t_0u_it_0u_j = \left\{\begin{array}{rl}
			u_{j+1}t_0u_it_0  & \mif 2 \le i \le j \\
			u_jt_0u_{i-1}t_0 & \mif 2 \le j < i 
			\end{array}\right.} \\[7pt]
	\end{array}$$
\end{thm}

\prf{}{
	See \cite{preprint4} for the proof as well as a discussion of the advantages and
	disadvantages over Grandis's presentation.
}

\rmk{ 
	Since, as mentioned in an earlier remark, all relations in the right
	column follow from the relations in the left column, all subsequent
	references to the statement of Theorem \ref{QuasiIds} will be understood
	as referring to relations of the left column only. 
}

\smallskip

Here are some other useful operators in $\fin^{op}$.

\defn{\label{replacementop}
	In the statement of Theorem \ref{QuasiIds}, note the overlapping conditions
	in the identities for $u_id_j$.  Indeed
	the equations
		$$u_id_i = d_{i+1}u_i = d_iu_{i+1} =: r_i$$
	hold for all $1 \le i \le n$.  We refer to the $r_i$ as {\it replacement operators}.
	They may also be defined directly as functions in $\fin$ as follows.
		$$r_i = r_i^{(n)}: [n] \lra [n] ~ \mfor n \ge 1 \mand 1 \le i \le n$$
		$$r_i(k) = \left\{\begin{array}{cl}
					0 & \mif k = 0 \mathrm{~or~} i \\
					k & \mathrm{~otherwise~} 
				\end{array}\right.$$
}

\begin{prp}
	For each $n \ge 1$, the replacement operators
		$$r_i: [n] \lra [n] \mfor 1\le i \le n$$
	constitute a family of mutually commuting idempotents in $\fin^{op}$.
		\begin{align*}
			r_i^2 &= r_i \\
			r_ir_j &= r_jr_i
		\end{align*}
\end{prp}
\prf{Proof.}{
	This is most easily verified using the formula for $r_i$ as a
	function in $\fin$ given in Definition \ref{replacementop}.
	Alternatively, it is a fun exercise to prove the assertion using the 
	algebraic identities of Theorem \ref{altpresoffin}.
}

\bigskip

%% file: Moorecomplex.tex

\label{Moorecomplex}
\label{SymmMooreCompActual}

In this section we briefly recall the construction, due to John C. Moore
\cite{JCMoore}, that takes a simplicial group $G$ and produces from it
its {\it Moore complex} $M(G)$. In this section we also describe extra
structure on the Moore complex of a symmetric-simplicial group in the
form of actions of permutation operators in $\fin^{op}$.

\defn{
	By the {\it Moore complex}
	$M = M(G)$ of the simplicial group $G$ we shall mean the chain complex
	consisting of the nonabelian groups
		$$M_n = M_n(G) := \setof{g\in G_n}{d_i(g) = 0 \mfor i \neq 0}$$
	and boundary operators
		$$d := d_0 : M_n \lra M_{n-1}.$$
	The group $M_n(G)$ will be called the group of {\it Moore $n$-chains} in $G$.
	Also for each $n\ge0$ the group of {\it Moore $n$-cycles} is defined as follows.
		\begin{align*}
			Z_0 &= Z_0(G) := M_0(G)\\
			Z_n &= Z_n(G) := \bigcap_{0\leq i \leq n} 
				~\mathrm{ker}~(d_i: G_n\rightarrow G_{n-1}) \text{~for $n\ge1$}
		\end{align*}
	Finally for each $n\ge0$ the group of {\it Moore $n$-boundaries} is
	defined as follows.
		\begin{align*}
			B_n &= B_n(G) := \text{Image}(d_0:M_{n+1}\lra G_n)
		\end{align*}
}

\rmk{
	Expositions about Moore
	complexes are to be found in \cite{May} and \cite{GoerssJardine}.
	Here we cite the following facts for the reader's edification.
	\begin{itemize}
		\item{$M(G)$ is a chain complex in the sense that $d_0^2:[n+2]\ra[n]$ is
			the trivial homomorphism for all $n\ge 0$.}
		\item{The group $B_n$ is normal in $Z_n$ for all $n\ge0$.}
		\item{The homology groups $Z_n/B_n$ of $M(G)$ are naturally isomorphic to the simplicial homotopy
			groups of $G$. }
	\end{itemize}
}

\bigskip

	We turn to consider the notion of Moore complex for symmetric-simplicial
	groups. Considering $G$ as a simplicial group (i.e., by restricting the
	action of $\fin^{op}$ to its subcategory $\ord^{op}$), one has the groups
	$M_n(G)$ and $Z_n(G)$ as defined above. Additionally, for all $n\ge0$, $G_n$ admits an action
	of $\sym_n$, the symmetric group on the set $[n] = \{ 0,1,\ldots,n\}$.
	Let 
		$$\sym^\p_n$$
	denote the subgroup of $\sym_n$ consisting of the permutations fixing $0$.

\begin{prp} 
	Let $G$ be a symmetric-simplicial group.
	The following hold with regard to the action of $\sym_n$ on $G_n$. 
	\begin{enumerate}
		\item{The subgroup $M_n(G)$ is invariant under the action of $\sym_n^\p$.}
		\item{The subgroup $Z_n(G)$ and therefore also $B_n(G)$ is invariant under the action of $\sym_n$.}
	\end{enumerate}
\end{prp}
\prf{Proof.}{
	For the first assertion, let $g$ belong to $M_n(G)$ so that all faces of $g$ except $d_0$ are trivial.
	Since $\sym_n^\p$ is generated by the transpositions $t_j$ for $j > 0$,
	it suffices to prove the following claim.
		$$d_it_j(g) = 0 \mforall i,j>0$$  
	An examination of the symmetric-simplicial identities (section \ref{SymmMooreComp}) 
	reveals that
		$$d_it_j = t_{j^\p}d_{i^\p} \mathrm{~or~} d_{i^\p}$$
	for some $i^\p$ and $j^\p$ where $i^\p \neq 0$, and from this the first assertion follows.

	The argument for the second assertion is similar. Taking $g\in Z_n(G)$, it suffices to prove
	the following.
		$$d_it_j(g) = 0 \mforall i,j \ge 0$$  
	This time, for {\it any} $i,j$, one has
		$$d_it_j = t_{j^\p}d_{i^\p} \mathrm{~or~} d_{i^\p}$$
	for some $i^\p$ and $j^\p$, and hence the second assertion follows.
}

\defn{\label{symmMoorecomp}
	Let $G$ be a symmetric-simplicial group.  By the {\it Moore complex of $G$}
	we shall mean the Moore complex $M(G)$ of the underlying simplicial group of $G$
	together with the actions of $\sym_n$ on $Z_n$ and $\sym^\p_n$ on $M_n$ and for all $n$.
}

The following definition abstracts the properties of symmetric Moore complexes.
Although we shall not need it here, we include it in order to highlight the fact that
the structures arising in this way are not the same as the notion of symmetric chain complex
encountered in the literature.

\defn{\label{abstsymmMoorecomp}
	Define a {\it symmetric chain complex} to be a chain complex of
	(nonabelian) groups
		$$\xymatrix{M_0 & M_1 \ar[l]_d & M_2 \ar[l]_d & \ldots \ar[l]}$$
	together with coextensive
	actions of $\sym_n$ on $Z_n := \mathrm{ker}(d:M_n\ra M_{n-1})$ 
	and $\sym^\p_n$ on $M_n$ for all $n$.  These data are required to satisfy 
	the following condition, which makes use of the generators $t_0,\ldots, t_{n-1}$ 
	of $\sym_n$ and in which $\sym^\p_n$ is identified 
	with the subgroup generated by $t_1,\ldots, t_{n-1}$.
	\begin{itemize}
		\item{\it For all $m\in M_n$ and all $i$ with $1 \le i \le n-1$, 
			$d(t_im) = t_{i-1}(dm)$.}
	\end{itemize}
}

\rmk{	The above condition is derived from the symmetric-simplicial identity
		$$d_0t_i = t_{i-1}d_0$$
	holding for $i\ge1$.  For $i=0$ one has the identity
		$$d_0t_0 = d_1$$
	which corresponds to the fact that in a symmetric chain complex, one always has
		$d(t_0m) = 0$
	for $m\in Z_n$ because $t_0$ preserves $Z_n$. 
}

%% file: DoldKanDecomps.tex

\label{DoldKanDecomps}

In this section, we generalize a result of
\cite{CarrascoCegarra} to the effect that the $n$th term $G_n$ of the
simplicial group $G$ is an internal \mbox{$2^n$}-SDP of
certain of its subgroups, each isomorphic to some term $M_j$ of the
Moore complex $M(G)$. These subgroups are either the normal subgroup
$M_n \lhd G_n$ or a copy of $M_{n^\p}$ for $n^\p < n$ embedded in $G_n$ as a 
subgroup of degenerate simplices via an iterated degeneracy operation as follows.
	$$\xymatrix{ M_{n-k} \subseteq G_{n-k} \ar[rr]^-{s_{i_k}s_{i_{k-1}}\ldots s_{i_1}} && G_n }$$

For $n\ge 3$, there is more than one order in which these subgroups can
be arranged to yield an SDP decomposition of $G_n$, and we give a
characterization of a large family of such orders, which includes (up to
choice of convention for the Moore complex) the total order discovered
in \cite{CarrascoCegarra}. Although we believe all such orders are
characterized in this fashion, we leave this question open (see Remark
\ref{dumbcharremark}).

\bigskip

Assume given a group $G$ and subgroups $H_1,\ldots,H_r$.

\defn{ \label{iSDP}
	The group $G$ is said to be an {\it internal $r$-semidirect product
	\emph{(briefly an {\it SDP})}~of the subgroups $H_i$} if the
	following two conditions hold.
	\begin{enumerate}
		\item{The set $H_1H_2\ldots H_i$ is a normal subgroup of $G$ for all $i$.}
		\item{Every $g\in G$ can be factored uniquely as a product 
			$$g = h_1h_2\ldots h_r$$ with $h_i\in H_i$ for all $i$.}
	\end{enumerate}
	We follow \cite{CarrascoCegarra} in using the notation
		$$G = H_1\rtimes\ldots\rtimes H_r$$
	if the above conditions hold.
}

\rmk{	The order in which the subgroups $H_i$ appear is an essential part of
	the definition of SDP. It can happen that $G$ is an SDP of the $H_i$
	when they arranged in certain orders, but not in others. At one extreme,
	there may be a unique such order. At the other extreme, it is
	immediate that $G$ is a direct product of the
	$H_i$ if and only if $G$ is an SDP of the $H_i$ arranged in each
	possible order.
}

\rmk{
	See \cite{preprint1} for a proof of the equivalence of Definition \ref{iSDP}
	with the original definition given in \cite{CarrascoCegarra}, as well as
	a further discussion of the properties of semidirect products and a
	corresponding higher external semidirect product construction. 
}

\bigskip

Let $G$ be a fixed simplicial group. Following \cite{CarrascoCegarra},
{\it additive notation is used for the group law of each term $G_n$,
although these groups are generally nonabelian.} Recall that $G$
possesses a Moore complex $M(G)$ whose terms will be denoted 
	$$M_n = M_n(G).$$

The following language and notations will be used. A {\it multi-index
$\alpha$ of length $k$ and dimension $\le n$} is a strictly increasing
sequence of indices
	$$\alpha = \big\{i_1 < \ldots < i_k \big\}$$
satisfying $0 \leq i_p \leq n-1$ for all $p$. The {\it length} of
$\alpha$ is the number $k$ of indices and is denoted by
$\vert\alpha\vert := k$. The set of all multi-indices of dimension $\le
n$ is denoted as follows.
	$$I^{(n)} := \big\{~\text{multi-indices $\alpha$ of dimension $\le n$}~\big\}$$

The {\it degeneracy operator} $s_\alpha$ corresponding to $\alpha$ is
the composition of elementary degeneracy operators
	$$s_\alpha : G_{n-\vert\alpha\vert} \longrightarrow G_n$$
	$$s_\alpha := s_{i_k}s_{i_{k-1}}\ldots s_{i_1}$$
and there is a corresponding subgroup of degenerate simplices
	$$H_\alpha := s_\alpha (M_{n-\vert\alpha\vert}) \subseteq G_n$$
for each $\vert\alpha\vert\ge 1$. Note this subgroup is an embedded copy
of $M_{n-\vert\alpha\vert}$, since the degeneracy operator $s_\alpha$ is
injective on account of the simplicial identities $d_is_i = d_{i+1}s_i =
\id$. For $\vert\alpha\vert=0$ there is only the \emph{empty multi-index
$\alpha = \varnothing$} and the corresponding subgroup
	$$H_\varnothing := M_n \lhd G_n$$
and its simplices (except the identity element) are all seen to be
nondegenerate by noting that any nontrivial simplex of $M_n$ has at most
one face not equal to the identity, whereas a nontrivial degenerate simplex always
has at least two equal nontrivial faces ($y = s_ix$ has the two equal
faces $d_iy = d_{i+1}y = x$).

\bigskip

Let us consider the question of how $G_n$ may decompose as an internal SDP of
the groups $H_\alpha$, that is, in which ways the $H_\alpha$ may be
totally ordered yielding such a decomposition. For example, if $G_n$ is
abelian, the $H_\alpha$ may be arranged in any order. More generally,
under various special conditions on the simplicial group $G$, there may
be especial flexibility in ordering the $H_\alpha$ to obtain
SDP decompositions of $G_n$. Here is a more specific question:
which total orders on the $H_\alpha$ give rise to SDP decompositions
of $G_n$ for all simplicial groups $G$, regardless of any special
characteristics of $G$? Such total orders will here be called {\it Dold-Kan total
orders}, and the corresponding universal decompositions {\it Dold-Kan
decompositions}.


Here we give an answer to this question as follows. There is a
partial order on the set $I^{(n)}$ of multi-indices of dimension $\le
n$, here called the {\it length-product partial order} for lack of a
better name, such that any total order extending the length-product
partial order will be a Dold-Kan total order. The proof is given in this
and the next section, and three examples of Dold-Kan total orders,
including the original one appearing in \cite{CarrascoCegarra}, are
given in the section after that. 

\rmk{\label{dumbcharremark}
	The question remains whether this includes {\it all} Dold-Kan total
	orders, that is, whether the Dold-Kan total orders are indeed
	characterized as those extending the length-product partial order. Although we
	believe the answer to this question is affirmative, we leave it open.
}

\defn{	For two multi-indices 
		$$\alpha = \big\{ i_1 < \ldots < i_k \big\} ~~~~~~~~
		\beta = \big\{ j_1 < \ldots < j_l \big\},$$
	both of degree $n$, define the {\it length-product} partial order relation 
		$$\alpha \prce \beta$$ 
	to mean	first of all that the condition $k \leq l$ holds, that is, 
	$$\begin{array}{rl}
		0. & \vert\alpha\vert \leq \vert\beta\vert,
	\end{array}$$
	and then that the following $k$ conditions are also satisfied:
	$$\begin{array}{rl}
		1. & i_1 \leq j_{l - k +1}; \\
		2. & i_2 \leq j_{l - k + 2}; \\
		& ~~~~\vdots \\
		k. & i_k \leq j_l
	\end{array}$$
	or equivalently that
		$$i_{k-p} \leq j_{l-p} \mforall 0\le p\le k-1.$$
The notation $\alpha \prec \beta$ indicates 
	the conjunction of $\alpha \prce \beta$ and $\alpha \neq \beta$.
}

Note the empty multi-index $\varnothing$ is the absolute minimum of $I^{(n)}$ under
the length-product partial order, and the full multi-index $\big\{ 0 < 1 <
\ldots < n-1 \big\}$ is the absolute maximum of $I^{(n)}$.

\bigskip

In the rest of this section, unless otherwise stated, $\alpha$ will
denote a fixed total order on the multi-indices of dimension $\le n$,
thought of as a bijection as follows.
	$$\alpha(\cdot): \big\{ 0, 1, \ldots, 2^n-1 \big\} \longrightarrow I^{(n)}$$
	$$k \mapsto \alpha(k)$$
On $\alpha$ the sole requirement is made that it extend the length-product partial order.  This
is equivalent to the requirement that $\alpha$ be {\it order-reflecting}, that is,
the following implication holds generally.
	$$\alpha(k) \prce \alpha(l) \implies k \leq l$$
In particular, $\alpha(0)$ is the empty multi-index $\varnothing$ and $\alpha(2^n-1)$ is 
the full multi-index.

\bigskip

The following theorem is the goal of the section.  It is convenient to
use the notation $H_k := H_{\alpha(k)}$.

\begin{thm} \label{maindecomp}
	There is an internal SDP decomposition as follows.
		$$G_n = H_0 \rtimes H_1 \rtimes \ldots \rtimes H_{2^n-1}$$
\end{thm}

The proof appears at the end of this section.  

\bigskip

In the meantime,
here are some tools for use in the proof. For an arbitrary
multi-index $\alpha$, the {\it face operator} corresponding in a dual
manner to $s_\alpha$ is
	$$d_\alpha^+ : G_n \longrightarrow G_{n-\vert\alpha\vert}$$
	$$d_\alpha^+ := d_{i_1}^+d_{i_2}^+\ldots d_{i_k}^+$$
where the following convenient notation is used.
	$$d_i^+ := d_{i+1}$$
There is also a corresponding {\it projection operator} defined as follows.
	$$\pi_\alpha: G_n \longrightarrow G_n$$
	$$\pi_\alpha(g) := s_\alpha d_{\alpha}^+(g)$$
Similarly as in the statement of Theorem \ref{maindecomp}, 
the notation $\pi_k := \pi_{\alpha(k)}$ is used.
Note that, by repeated application of the simplicial
identities $d_{i+1}s_i = \id$, one obtains the identity
	\begin{equation} \label{supcancid}
		d_\alpha^+s^{\phantom+}_\alpha = \id
	\end{equation}
and consequently $\pi_k$ is indeed a projection in the 
sense that the following holds.
	$$\pi_k^2 = \pi_k$$


Here is an algorithm that will be shown to produce, for an
arbitrary element $g\in G_n$, a decomposition of the kind indicated in
Theorem \ref{maindecomp}. Define a sequence of elements $g_i \in G_n$,
starting with $g_{2^n-1} := g$ and proceeding recursively in reverse
order by the formula 
	\begin{equation} \label{recurse} g_{k-1} := g_k - \pi_k(g_k)
	\end{equation}
and ending with $g_0$.
By induction one has
	\begin{align*} g_k & = g_{k-1} + \pi_k(g_k) \\
		       	 & = g_{k-2} + \pi_{k-1}(g_{k-1}) + \pi_k(g_k) \\
			 & \hspace{6pt} \vdots \\
			 & = g_0 + \pi_1(g_1) + \pi_2(g_2) + \ldots + \pi_k(g_k)
	\end{align*}
which, following \cite{CarrascoCegarra}, we also denote
	\begin{equation*}
		g_k = \sum_{i = 0}^k \pi_i(g_i)
	\end{equation*}
keeping in mind that the addition here is not commutative in general, so
that it is essential to order the summands with indices increasing from
left to right. Thus one has the following decomposition for any $g\in
G_n$.
	\begin{equation} \label{factorize}
		g = g_{2^n-1} = \sum_{i = 0}^{2^n-1} \pi_i(g_i)
	\end{equation}
Equation \eqref{factorize} will be referred to as the {\it $\alpha$-decomposition}
of $g$.
\bigskip

The proof of the SDP decomposition of Theorem \ref{maindecomp} above requires
certain generalized versions of the simplicial identities as well as an
understanding of how the length-product partial order arises from them.
This is accomplished by the following theorem, which is proved in detail
in the next section.

\bigskip
\begin{thm} \label{ids}
	The following facts hold for arbitrary multi-indices 
		$$\alpha = \big\{ i_1 < \ldots < i_k \big\} ~~~~~~~~
		\beta = \big\{ j_1 < \ldots < j_l \big\}$$
		of dimension $\le n$.
	\begin{enumerate}
		\item{If $d_{\alpha}^+ s^{\phantom+}_\beta = s^{\phantom+}_{\beta^\prime}$ for some multi-index $\beta^\prime$, 
			then $\alpha \prce \beta$.}
		\item{If $\alpha \nprce \beta$ then $d_{\alpha}^+ s^{\phantom+}_\beta 
				= d_{\alpha^\prime}^+s^{\phantom+}_{\beta^\prime}d^{\phantom+}_i$
				for some $\alpha^\prime,\beta^\prime$ and some $i \neq 0$.}
		\item{For $i \neq 0$, one has $d_id_{\beta}^+ = d_{\beta^\prime}^+$
			for some multi-index $\beta^\prime \succ \beta$.}
	\end{enumerate}
\end{thm}

\bigskip

The following lemmas refer to the $\alpha$-decomposition of an element $g\in G$
given by \eqref{recurse}, \eqref{factorize}.

\begin{lem} \label{kernel}
	The following holds for any $k$ with $0 \leq k \leq 2^n-1$.
		$$g_k \in 
			\bigcap_{i = k+1}^{2^n-1}\mathrm{ker}~d_{\alpha(i)}^+ = 
			\bigcap_{i = k+1}^{2^n-1}\mathrm{ker}~\pi_{i}$$
\end{lem}
\prf{Proof.}{
	Proceed by induction (in reverse order).  To start, note that
	the case $k=2^n-1$ is vacuously true.
\newline\indent	Now assume the statement holds for some $k$ with $0 < k \leq 2^n-1$.
	One checks as follows for any $i \ge k$ that $d_{\alpha(i)}^+(g_{k-1})$ is trivial.
	For $i = k$, one has
		\begin{align*}
			d_{\alpha(k)}^+(g_{k-1}) & = d_{\alpha(k)}^+(g_k) - 
				d_{\alpha(k)}^+\Big(\pi_k(g_k)\Big) \tag{By \eqref{recurse}}\\
			&=  d_{\alpha(k)}^+(g_k) - 
				d_{\alpha(k)}^+\Big(s^{\phantom+}_{\alpha(k)}d_{\alpha(k)}^+(g_k)\Big) \\
			&= d_{\alpha(k)}^+(g_k) - d_{\alpha(k)}^+(g_k) \tag{By \eqref{supcancid}}\\ 
			&= 0.
		\end{align*}
	For $i > k$ one starts out similarly, obtaining
		$$ d_{\alpha(i)}^+(g_{k-1}) = d_{\alpha(i)}^+(g_k) - 
			d_{\alpha(i)}^+\Big(s^{\phantom+}_{\alpha(k)}d_{\alpha(k)}^+(g_k)\Big)$$
	of which the first term is trivial by the inductive hypothesis.  To see
	that the second	term is also trivial, apply part 2 of Theorem \ref{ids} to find
		$$d_{\alpha(i)}^+s^{\phantom+}_{\alpha(k)}d_{\alpha(k)}^+ = 
			d_{\alpha^\p}^+ s^{\phantom+}_{\beta^\p} d_i^{\phantom+} d_{\alpha(k)}^+$$
	for some multi-indices $\beta$ and $\gamma$ and some index $i \neq 0$.  Hence by part 3 of
	Theorem \ref{ids} one has 
		$$d_{\alpha(i)}^+s^{\phantom+}_{\alpha(k)}d_{\alpha(k)}^+ 
			= d_{\alpha^\p}^+ s^{\phantom+}_{\beta^\p} d_{\alpha(k^\prime)}^+$$
	where $\alpha(k^\p)\succ\alpha(k)$ and hence $k^\p > k$.  
	Therefore, again by the inductive hypothesis, this operator annihilates $g_k$
	as claimed.
}


\begin{lem} \label{pikernel}
	The following holds for any $k$ with $0 \leq k \leq 2^n-1$.
		$$\pi_k(g_k) \in H_k$$
\end{lem}
\prf{Proof.}{
	The assertion of the lemma is that
		$$s^{\phantom+}_{\alpha(k)}d_{\alpha(k)}^+(g_k) \in s_{\alpha(k)}M_{n-\abs{\alpha(k)}}$$
	or equivalently
		$$d_{\alpha(k)}^+(g_k) \in M_{n-\abs{\alpha(k)}}$$
	that is, it suffices to show that $d^{\phantom+}_i d_{\alpha(k)}^+(g_k)$ trivial for any $i\neq 0$.
	By part 3 of Theorem \ref{ids}, one has
		$$d^{\phantom+}_i d_{\alpha(k)}^+ = d_{\alpha(k^\prime)}^+$$
	for some $\alpha(k^\p)\succ\alpha(k)$ and hence $k^\p > k$, and so the operator
	$d_{\alpha(k^\prime)}^+$ annihilates $g_k$ by Lemma \ref{kernel}.  Hence 
	$d_i$ annihilates $d_{\alpha(k)}^+(g_k)$, as required.
}

\begin{lem}\label{ids8}
	For any $k^\prime > k$, the operator $\pi_{k^\prime}$ annihilates $H_k$.
\end{lem}
\prf{Proof.}{
	It suffices to show that $d_{\alpha(k^\prime)}^+$ annihilates
	$s^{\phantom+}_{\alpha(k)}M^{\phantom+}_{n-k}$. Since
		$$k^\p > k \implies \alpha(k^\prime) \nprce \alpha(k)$$
	one has by part 2 of Theorem \ref{ids}
		$$d_{\alpha(k^\prime)}^+ s^{\phantom+}_{\alpha(k)}M_{n-k} = 
			d_{\alpha^\p}^+ s^{\phantom+}_{\beta^\p} d^{\phantom+}_iM_{n-k}$$
	for some multi-indices $\beta$ and $\gamma$ and some index $i \neq 0$.  But
		$$d_iM_{n-k} = \{0\}$$
	by definition of the group $M_{n-k}$, and the claim follows.
}

\bigskip

\prf{Proof of Theorem \ref{maindecomp}.}{
	One must verify the two conditions of Definition \ref{iSDP}.  
	First it is claimed that, for any $g\in G_n$, its
	$\alpha$-decomposition \eqref{factorize} is the unique factorization of $g$ of
	the form
		$$g = h_0 + h_1 + \ldots + h_{2^n-1}$$
	with $h_i \in H_i$ for all $i$.  That the $\alpha$-decomposition of $g$
	is such a factorization is just Lemma \ref{pikernel}.
	For uniqueness, let $g\in G_n$ have two factorizations
		$$g = \sum_{i=0}^{2^n-1}h_i = \sum_{i=0}^{2^n-1}h_i^\prime$$
	in which $h_i$ and $h_i^\prime$ belong to $H_i$ for all $i$.
	By Lemma \ref{ids8}, applying the projection
	$\pi_{2^n-1}$ to both sides yields
		$$h_{2^n-1} = h_{2^n-1}^\prime$$
	and cancelling these terms on the right, one applies the next 
	projection $\pi_{2^n-2}$ to get the next right-most terms equal.
	Inductively it follows for all $i$ that
		$$h_i = h_i^\prime.$$
	
	To verify the other condition of the definition of internal $r$-SDP, 
	namely that $H_0+H_1+\ldots +H_k$ is a normal subgroup of $G$ for
	each $k$, one must prove the following equality, due in its original form
	to \cite{CarrascoCegarra}.
		$$H_0+H_1+\ldots +H_k = \bigcap_{i = k+1}^{2^n-1}{\text{ker~}}\pi_{i}$$
	From the uniqueness statement just proved,
	one deduces the following for the $\alpha$-decomposition
	of an element $g\in G$ that belongs to the left-hand side.
		\begin{align*}
			g=h_0 + \ldots + h_k &\iff \pi_{k^\p}(g_{k^\p}) = 0 ~\mforall~ k^\p > k \\
			&\iff g_{k^\p} = g ~\mforall~ k^\p \ge k 
		\end{align*}
	A simple induction argument shows that the last two statements together are equivalent
	to 
		$$\pi_{k^\p}(g) = 0 ~\mforall~ k^\p > k$$
	which says that $g$ belongs to the right-hand side above, as claimed.
}

%% file: ProofLemmaIdsNew.tex

	\label{ProofLemmaIds}	\label{ProofLemmaIdsNew}

The purpose of this section is to prove in detail Theorem \ref{ids},
used in the last section. For convenience, we recall one of the
simplicial identities in a slightly modified form. \\

\begin{thing}{A Simplicial Identity.}
	For $n \ge 0$ and $0 \leq i,j \leq n$, the simplicial operator
	$$d_i^+ s_j : G_n \longrightarrow G_n$$
	can be rewritten as follows.
	$$\boxed{\newline ~~d_i^+ s_j = \left\{ \begin{array}{cl} 	s_{j-1}d_i^+ & \text{if}~~i < j-1 \\
						\text{id}  &\text{if}~~i = j-1, j \\
						s_j d_{i-1}^+&\text{if}~~i > j
				\end{array}\right.}$$ 
	Note the formulas have the same appearance as the corresponding usual
	identity except for a shift in the conditions, that is, $i$ has been
	replaced by $i+1$. The cases corresponding to $d_0s_j$ are absent---although 
	of course still true---reflecting the special role of $d_0$ as
	boundary operator in our chosen convention of Moore complex.
\end{thing} 

\bigskip

The following convenient notations will be used.
	$$s_j^- := s_{j-1}$$
	$$s_\beta^- := s_{j_l}^-\ldots s_{j_1}^-$$

\bigskip

The heart of the proof of Theorem \ref{ids} is given in the following proposition.
	
\begin{prp} \label{idslemma} 
	Let $\beta = \big\{ j_1 < \ldots < j_l \big\}$ 
	be a multi-index of dimension $\le n$ and let $i$ be an index satsifying $0 \le i \le n-1$.
	Consider the following two alternatives, which are mutually exclusive and
	cover all possibilities for $\beta$ and $i$.  (The nonsensical statements
	$j_{l+1}-1 > i$ and $i > j_0$ are regarded as true for any $i$.)
	\begin{enumerate}
		\item{For some (unique)	subscript $0 \le q \le l$, one has $j_{q+1}-1 > i > j_q$.}
		\item{For some (unique) subscript $1 \le q \le l$, one has
			$j_{q+1}-1 > i$ and $i = j_q$ or $j_q-1$.}
	\end{enumerate}	
	In the first case, the identity 
			$$d_i^+s_\beta = s_{j_l}^-s_{j_{l-1}}^-\ldots 
				s_{j_{q+1}}^-s_{j_q}^{\phantom{-}}\ldots s_{j_1}^{\phantom{-}} d_{i-q}^+$$
	holds, and the rightmost factor $d_{i-q}^+$ is different from $d_0$.
	Then we say that $d_i^+$ has ``slipped past'' $s_\beta$. 
	In the second case, the identity 
			$$d_i^+s_\beta = s_{j_l}^-s_{j_{l-1}}^-\ldots 
				s_{j_{q+1}}^-s_{j_{q-1}}^{\phantom{-}}\ldots s_{j_1}^{\phantom{-}}$$
	holds, and we say that $d_i^+$ is ``absorbed by'' $s_\beta$. 

\end{prp}


The following lemma, whose statement uses the same
notations, is necessary for the proof.

\begin{lem} \label{idslemmalemma}
	The following special cases of Proposition \ref{idslemma} hold.
	\begin{enumerate}
		\item{If $i>j_l$ then $d_i^+s_\beta^{\phantom{+}} = s_\beta^{\phantom{+}} d_{i-\abs{\beta}}^+$.}
		\item{If $i < j_1 - 1$ then $d_i^+s_\beta^{\phantom{+}} = s_\beta^- d_i^+$.}
	\end{enumerate}
	In each case, the face operator at the right (that is, $d_{i-\abs{\beta}}^+$ 
	or $d_i^+$) is not $d_0$.
\end{lem}
\prf{Proof.}{
	For the first assertion, note that 
	$i > j_l$ implies on account of the strictly 
	increasing nature of the indices of $\beta$
		$$i - p > j_l - p \ge j_{l-p}$$
	for any $p$ with $0 \leq p \leq l-1$.  Hence the simplicial identity above
	may be applied repeatedly and one obtains thus
	\begin{equation*}
	\begin{split}
		d_i^+s_\beta & = d_i^+s_{j_l}s_{j_{l-1}}\ldots s_{j_1} \\
		& = s_{j_l}d_{i-1}^+s_{j_{l-1}}\ldots s_{j_1} \\
		& = s_{j_l}s_{j_{l-1}}d_{i-2}^+s_{j_{l-2}}\ldots s_{j_1} \\
		& \hspace{6pt} \vdots \\
		& = s_{j_l}s_{j_{l-1}}\ldots s_{j_1} d_{i-l}^+ \\
		& = s_\beta d_{i-\vert\beta\vert}^+
	\end{split}
	\end{equation*}
	as claimed.  To see that $d_{i-\vert\beta\vert}^+$ is not $d_0$, note
	that because the $j_q$ are strictly increasing one has $j_q \ge q-1$ for each $q$.
	From the hypothesis $i > j_l$ it then follows that $i \ge \abs{\beta}$
	and so indeed $d_{i-\abs{\beta}}^+$ cannot be $d_0$.

\bigskip	

	For the second assertion, note it follows from the increasing
	nature of the indices of $\beta$ and the hypothesis $i < j_1 -1$ that
		$$i < j_q - 1$$
	for all $q$.  Then, again successively applying the simplicial identity, one obtains
	\begin{equation*}
	\begin{split}
		d_i^+s_\beta & = d_i^+s_{j_l}s_{j_{l-1}}\ldots s_{j_1} \\
		& = s_{j_l}^-d_i^+s_{j_{l-1}}\ldots s_{j_1} \\
		& = s_{j_l}^-s_{j_{l-1}}^-d_i^+s_{j_{l-2}}\ldots s_{j_1} \\
		& \hspace{6pt} \vdots \\
		& = s_{j_l}^-s_{j_{l-1}}^-\ldots s_{j_1}^-d_i^+ \\
		& = s_\beta^- d_i^+
	\end{split}
	\end{equation*}
	as claimed.  From $i \ge 0$ it immediately follows that $d_i^+$ is not $d_0$.
%
}

\prf{Proof of Proposition \ref{idslemma}.}{
	If $\beta$ contains neither $i$ nor $i+1$, then the situation of part 1
	of the Proposition obtains.  In that case, there is a factorization
		$$s_\beta = s_{\beta^\p}s_{\beta^\pp}$$
	where $\beta^\p$ consists of those indices of $\beta$ greater than $i+1$ and
	$\beta^\pp$ consists of those indices of $\beta$ less than $i$.  
	Applying Lemma \ref{idslemmalemma}, one obtains
		\begin{align*} 
			d_i^+ s_\beta^{\phantom{+}} &= d_i^+s_{\beta^\p}^{\phantom{+}}s_{\beta^\pp}^{\phantom{+}} \\
			&= s_{\beta^\p}^- d_i^+s_{\beta^\pp}^{\phantom{+}} \tag{By Lemma \ref{idslemmalemma}, part 1}\\
			&= s_{\beta^\p}^- s_{\beta^\pp}^{\phantom{+}} d_{i-\abs{\beta^\pp}}^+ \tag{By Lemma \ref{idslemmalemma}, part 2}
		\end{align*}
	where $d_{i-\abs{\beta^\pp}}^+$ is not $d_0$.  This proves part 1 of the Proposition.

	Now assume $\beta$ contains one or both of $i$ and $i+1$, so that the
	situation of part 2 of the Proposition obtains.  Letting $p$ stand for the larger
	of $i$ and $i+1$ contained in $\beta$, there is a factorization
		$$s_\beta = s_{\beta^\p}s_p s_{\beta^\pp}$$
	where $\beta^\p$ consists of those indices of $\beta$ greater than $p+1$ and
	$\beta^\pp$ consists of those indices of $\beta$ less than $p$.  
	Again applying Lemma \ref{idslemmalemma}, one obtains
		\begin{align*}
			d_i^+ s_\beta &= d_i^+s_{\beta^\p}^{\phantom{+}}s_p^{\phantom{+}} s_{\beta^\pp}^{\phantom{+}} \\
			&= s_{\beta^\p}^- d_i^+ s_p^{\phantom{+}} s_{\beta^\pp}^{\phantom{+}} \tag{By Lemma \ref{idslemmalemma}, part 1}\\
			&= s_{\beta^\p}^-s_{\beta^\pp}^{\phantom{+}} \tag{By the simplicial identity}
		\end{align*}
	thus proving part 2 of the Proposition.
}

We turn to the proof of Theorem \ref{ids} from the previous section.
Recall that there the rank-product partial order $\prce$ was defined on
the multi-indices of dimension $\le n$.

\bigskip
\thing{Theorem \ref{ids}}{~~The following facts hold for arbitrary multi-indices 
		$$\alpha = \big\{ i_1 < \ldots < i_k \big\} ~~~~~~~~
		\beta = \big\{ j_1 < \ldots < j_l \big\}$$
		of dimension $\le n$.
	\begin{enumerate}
		\item{If $d_{\alpha}^+ s^{\phantom+}_\beta = s^{\phantom+}_{\beta^\prime}$ for some multi-index $\beta^\prime$, 
			then $\alpha \prce \beta$.}
		\item{If $\alpha \nprce \beta$ then $d_{\alpha}^+ s^{\phantom+}_\beta 
				= d_{\alpha^\prime}^+s^{\phantom+}_{\beta^\prime}d^{\phantom+}_i$
				for some $\alpha^\prime,\beta^\prime$ and some $i \neq 0$.}
		\item{For $i \neq 0$, one has $d_id_{\beta}^+ = d_{\beta^\prime}^+$
			for some multi-index $\beta^\prime \succ \beta$.}
	\end{enumerate}
}\bigskip

\prf{Proof of Theorem \ref{ids}, part 1.}{
	The strategy is to work with the left hand side of 
		$$d_{\alpha}^+ s_\beta = s_{\beta^\prime}$$
	using the simplicial identities to push the elementary factors of
	$d_\alpha^+$ one at a time across the sequence of factors of $s_\beta$
	and to observe in the process that the various conditions constituting
	the assertion $\alpha \prce \beta$ hold.

	First observe that, according to the dichotomy of Proposition \ref{idslemma},
	as each factor of $d_\alpha^+$ is pushed through $s_\beta$, it will
	either cancel a factor of $s_\beta$ or it will slip past $s_\beta$. In
	the present case, it is not possible for a factor to slip past, because
	any product of elementary operators starting with a face operator $d_i$
	on the right corresponds to a monotonic function in \ord~failing to have
	$i$ as a value, whereas $s_{\beta^\prime}$ corresponds to a surjective
	function in \ord. One concludes that each factor of $d_\alpha^+$ cancels
	some factor of $s_\beta$, and consequently
		$$\vert\alpha\vert \leq \vert\beta\vert.$$

	To verify the remaining conditions, now begin by pushing $d_{i_k}^+$
	across $s_\beta$. Under	Proposition \ref{idslemma}, 
	one may say that there is a unique subscript
	$q(k)$ with $0 \leq q(k) \leq l$ satisfying
		$$j_{q(k)+1} - 1 > i_k$$
		\begin{equation} \label{qkdef}
			i_k = j_{q(k)} \text{~~or~~} j_{q(k)}-1 
		\end{equation}
	so that $d_{i_k}^+$ cancels with $s_{j_{q(k)}}$, and moreover the result is
		$$d_{i_k}^+s_\beta = s_{j_l}^-s_{j_{l-1}}^-\ldots 
				s_{j_{q(k)+1}}^-s_{j_{q(k)-1}}\ldots s_{j_1} =: s_{\beta(k)}.$$

	Now we seek to push the next face operator $d_{i_{k-1}}^+$ across
	$s_{\beta(k)}$.  Since $d_{i_k}^+$ slipped past $s_{j_l},
	\ldots, s_{j_{q(k)+1}}$, the next factor $d_{i_{k-1}}^+$ will also slip past $s_{j_l}^-,
	\ldots, s_{j_{q(k)+1}}^-$ on account of the strict inequality $i_{k-1} <
	i_k$. To be precise, one has the conditions
		$$i_{k-1} < i_k <  j_{q(k)+1}-1 < \ldots < j_l-1$$
		$$\implies i_{k-1} < j_{q(k)+1}-2 < \ldots < j_l-2$$
	enabling us to apply the simplicial identity to get
		$$d_{i_{k-1}}^+s_{\beta(k)} = s_{j_l}^{--}s_{j_{l-1}}^{--}\ldots 
				s_{j_{q(k)+1}}^{--}d_{i_{k-1}}^+s_{j_{q(k)-1}}\ldots s_{j_1}$$
	where $s^{--}_j := s_{j-2}$. As 
	$d_{i_{k-1}}^+$ is pushed further to the right, again there must be a
	unique subscript $q(k-1)$, evidently less than $q(k)$, such that
	$d_{i_{k-1}}^+$ cancels with $s_{j_{q(k-1)}}$.  Again by Proposition \ref{idslemma}, the result is
	\begin{multline*}
		d_{i_{k-1}}^+s_{\beta(k)} = s_{j_l}^{--}s_{j_{l-1}}^{--}\ldots s_{j_{q(k)+1}}^{--}\circ \\
		s_{j_{q(k)-1}}^-s_{j_{q(k)-2}}^-\ldots s_{j_{q(k-1)+1}}^-\circ \\
		s_{j_{q(k-1)-1}}s_{j_{q(k-1)-2}}\ldots s_{j_1} =: s_{\beta(k-1)}.
	\end{multline*}

	Continuing in this manner, one obtains a sequence of subscripts $q(p)$ with
		\begin{equation} \label{qmonofunc}
			l \ge q(k) > q(k-1) > \ldots > q(1) \ge 0
		\end{equation}
	and a sequence of degeneracy operators $s_{\beta(p)}$ such that 
		$$d_{i_p}^+ s_{\beta(p+1)} =:  s_{\beta(p)}$$
	and such that in
	pushing $d_{i_p}^+$ across $s_{\beta(p+1)}$, it cancels with
	$s_{j_{q(p)}}$ and hence does not affect the indices of the operators
	$s_{j_q}$ further to the right (that is, the indices $q < q(p)$).  Since the
	sequence of indices $j_q$ is strictly increasing, one has
		\begin{align*}
			i_p &\leq j_{q(p)} \tag{By \eqref{qkdef}}\\
			&= j_{q(k - k + p)} \\
			&\leq j_{q(k)-k+p} \tag{By \eqref{qmonofunc}} \\
			&\leq j_{l-k+p} \tag{since $q(k) \le l$}
		\end{align*}
	for all $p$ with $1 \leq p \leq k$, and thus it is verified that $\alpha \prce \beta$.
}

\prf{Proof of Theorem \ref{ids}, part 2.}{
	Recall that, by Proposition \ref{idslemma}, in pushing each factor of $d_\alpha^+$ across
	$s_\beta$, one of two possibilities can occur: either the factor cancels
	somewhere along the way, or it makes it all the way through to the
	right. If cancellation occurred for each factor, then by the just-proved
	part 1 of the Theorem, one would have $\alpha \prce \beta$. Since it was assumed that
	$\alpha \nprce \beta$, one of the factors must make it through.

	Using the notation from the proof of part 1, let us say that $d_{i_p}^+$ is the
	first factor to make it through.  Prior to this occurrence, one has
		\begin{align*}
			d_\alpha^+ s_\beta & = d_{i_1}^+\ldots\ldots d_{i_{k-2}}^+d_{i_{k-1}}^+d_{i_k}^+ s_\beta \\
			& = d_{i_1}^+\ldots\ldots d_{i_{k-2}}^+d_{i_{k-1}}^+ s_{\beta(k)} \\
			& = d_{i_1}^+\ldots\ldots d_{i_{k-2}}^+ s_{\beta(k-1)} \\
			& \hspace{6pt} \vdots \\
			&  = d_{i_1}^+ \ldots d_{i_{p-1}}^+d_{i_p}^+s_{\beta(p+1)}
		\end{align*}
	The hypotheses of Proposition \ref{idslemma}, part 1 
	apply to $d_{i_p}^+s_{\beta(p+1)}$, for otherwise there 
	would be cancellation of $d_{i_p}^+$.  Therefore one gets finally
		$$d_\alpha^+ s_\beta = d_{i_1}^+ \ldots d_{i_{p-1}}^+s_{\beta^\prime}d_i^+$$
	where $i \neq 0$.
}

\prf{Proof of Theorem \ref{ids}, part 3.}{
	For this final proof, we shed our previous labelling habits and
	index $d_\beta^+$ directly, writing
		$$d_\beta^+ = d_{j_1}d_{j_2}\ldots d_{j_l}$$
		$$1 \leq j_1 < \ldots < j_l \leq n$$
	and consider how the simplicial identities give rise to the multi-index $\beta^\p$
	in the equation $d_{\beta^\p}^+ = d_id_\beta^+$.  The relevant computation
		\begin{align*}
			d_id_\beta^+ & = d_{i+0}d_{j_1}d_{j_2}d_{j_3}\ldots\ldots\ldots d_{j_l} \\
			& = d_{j_1} d_{i+1}d_{j_2}d_{j_3}\ldots\ldots\ldots d_{j_l} \\
			& = d_{j_1}d_{j_2} d_{i+2}d_{j_3}\ldots\ldots\ldots d_{j_l} \\
			& \hspace{6pt} \vdots \\
			& = d_{j_1}\ldots d_{j_q} d_{i+q} d_{j_{q+1}} \ldots d_{j_l}
		\end{align*}
	is explained as follows.
	As long as 
		$$i+m \ge j_{m+1}$$ 
	holds, one may push $d_{i+m}$ (the avatar of $d_i$) to the right using the simplicial identity
		$$d_{i+m}d_{j_{m+1}} = d_{j_{m+1}}d_{i+m+1}$$
	until the first subscript $q$ is reached such that
		$$i+q < j_{q+1},$$
	at which point the indices are in increasing order from left to right.  
	Since the product on the right-hand-side evidently does not contain $d_0$ as a factor,
	there exists a multi-index $\beta^\p$ such that the final result $d_{\beta^\prime}^+$.
	Now reindex it directly as
		$$d_{\beta^\prime}^+ = d_{j_1^\prime} d_{j_2^\prime} \ldots d_{j_{l+1}^\prime}$$
		$$j_p^\prime = 	\left\{\begin{array}{cl}
					j_p & \text{~~if~~} 1 \leq p \leq q \\
					i+q & \text{~~if~~} p = q+1 \\
					j_{p-1} & \text{~~if~~} p \ge q+2.
				\end{array}\right.$$
	It is straightforward to verify $j_p \leq j_{1+p}^\prime$ for each $p$ with $1\leq p \leq l$,
	showing $\beta \prce \beta^\prime$.  
	Due to the fact that
		$$\vert\beta^\prime\vert = \vert\beta\vert + 1$$
	it must be that $\beta \neq \beta^\p$ and so one concludes $\beta \prec \beta^\prime$ as claimed.
}


%% file: BinaryTO.tex

	\label{BinaryTO}

In this section, we describe the Dold-Kan total order discovered by
Carrasco and Cegarr\`a \cite{CarrascoCegarra}. The difference in its
appearance here is due to our choice of convention for the Moore complex
(see section \ref{Moorecomplex}). Written this way, it is given by
binary representations of the natural numbers, and so it seems
reasonable to call it the {\it binary order}. 

\bigskip

\defn{ \label{binaryorder}
	The {\it binary total order} $\alpha(\cdot)$ is defined as follows.
	For a nonnegative integer $k$ with binary expansion
		$$k = \sum_{i=0}^{n-1} c_i2^i$$
	where $c_i$ is the $i$th binary digit of $k$, set
		$$\alpha(k) := \big\{~ i ~\big\vert~ c_i = 1 ~\big\}.$$
	For $n=4$ the total order appears as follows.  The indices are written
	as they appear in the subscripts of degeneracy operators, that is,
	in decreasing order.
	$$\begin{array}{cccccccc}
		&\varnothing & < & 0 & < & 1 & < & 10 \\
		 <& 2 & < & 20 & < & 21 & < & 210 \\
		 <& 3 & < & 30 & < & 31 & < & 310 \\
		 <& 32& < & 320& < & 321& < & 3210
	\end{array}$$
	In general, $\alpha^\p$ is strictly less than $\alpha$ in
	the binary order if and only if the inequality
		$$\sum_{j\in\alpha^\p}2^j < \sum_{j\in\alpha}2^j$$
	holds, that is, in reading the corresponding 
	binary expansions
	from left to right (greatest to least), if the $p$-th digit
	is the first place in which they differ, then the $p$-th digit of $\alpha^\p$ is a 0
	and for $\alpha$ it is 1 (i.e., $p$ belongs to $\alpha$ but not to $\alpha^\p$).
}

\rmk{
	A different, but equivalent, description is that $\alpha^\p$ is strictly
	less than $\alpha$ in the binary order if and only if the corresponding
	degeneracy operators are of the form
		\begin{align*}
			s_{\alpha^\p} &= s_{i_k}s_{i_{k-1}} \ldots s_{i_{q+1}} s_{i^\p_q} s_{i^\p_{q-1}} \ldots \\
			s_\alpha	 &= s_{i_k}s_{i_{k-1}} \ldots s_{i_{q+1}} s_{i_q} s_{i_{q-1}} \ldots
		\end{align*}
	where the index $i^\p_q$ is either strictly less than $i_q$ or
	nonexistent (we assume the labelling here to be such that either
	$\alpha^\p$ or $\alpha$ may have rank $k$ or less). That is, comparing
	the indices of $\alpha^\p$ and $\alpha$ in order from greatest to least,
	either $\alpha$ has the greater index in the first position $q$ in which
	they differ, or they coincide in each position up until the point that
	$\alpha^\p$ runs out of indices and $\alpha$ still has some left.
\newline\indent	
	Essentially for this reason, Carrasco and Cegarra \cite{CarrascoCegarra}
	call this order {\it lexicographic}, and it also appears this way in the
	indices of degeneracy operators, that is, when the indices are written
	in decreasing order, as was done above. According to our chosen
	conventions, however (namely, our habit of writing multi-indices as increasing from left to right
	and also our choice of convention for the Moore complex), binary order would be called {\it
	reverse-lexicographic}. We use the name {\it binary} partly in order to
	avoid confusion on this point. 
}

\rmk{\label{CCremarkableform}
	Carrasco and Cegarra discovered a remarkable explicit formula for the
	individual components $\pi_k(g_k)$ (see section \ref{DoldKanDecomps})
	holding when binary order is used, and the reader is referred to their paper
	\cite{CarrascoCegarra} for its proof. The purpose of this remark is to give one
	way of stating their formula using the present conventions. 
	The notations of the section \ref{DoldKanDecomps} will be freely used.
\newline\indent	
	Fix an index $m$ with $0 \le m \le 2^n-1$, and let 
		$$\alpha = \{i_1 < \ldots < i_k\}$$ 
	stand for $\alpha(m)$. Also write $\alpha^c$ for the multi-index which
	is the complement of $\alpha$ in $\{0,\ldots,n-1\}$ and write
		$$\alpha^c = \{j_1 < \ldots < j_l\}.$$  
	Define the homomorphisms
		$$q_i: G_n \lra G_n$$
		$$q_i(g) := s_i d_i^+(g)$$
	for $0\le i \le n-1$, and let
		$$q_\alpha := q_{i_k} \ldots q_{i_1}$$
	be the product of $q_i$ for $i\in\alpha$ in {\it decreasing} order from
	left to right. Also write $q_j^\perp$ for the crossed homomorphisms
		$$q_j^\perp: G_n \lra G_n$$
		$$q_j^\perp := 1 - q_j$$
		$$q_j^\perp(g) = g - s_jd_j^+(g)$$
	and write $q_{\alpha^c}^\perp$ for the product 
		$$q_{j_1}^\perp\ldots q_{j_l}^\perp$$
	of $q_j$ for $j\in\alpha^c$ in {\it increasing} order
	from left to right.
\newline\indent	
	Then the component $\pi_k(g_k)$ of $g\in G_n$ under the SDP
	decomposition of the previous section using binary order can be
	expressed explicitly in terms of $g$ by the following formula.
		$$\pi_k(g_k) = q_\alpha^{\phantom{\perp}} q_{\alpha^c}^\perp(g)$$
}

%% file: SymmDoldKan.tex

	\label{SymmDoldKan}

In this section, we show that, if $G$ is a symmetric-simplicial group,
there are new Dold-Kan-type decompositions available for it in addition
to the ones described in section \ref{DoldKanDecomps}. Once again the
$n$th term $G_n$ is a \mbox{$2^n$-fold} internal SDP of certain of its
subgroups, each isomorphic to some term $M_j$ of the Moore complex
$M(G)$. The main difference is that the subgroups $s_\alpha
M_{n-\abs{\alpha}}$ of degenerate simplices are replaced by subgroups
$u_\alpha M_{n-\abs{\alpha}}$ of quasidegenerate simplices.

As to the question of the ordering of these subgroups, many more
orderings are available than in the earlier case. In effect, the
length-product partial order is replaced by the partial order given by
inclusion, and requiring a total order to extend this partial order
places many fewer constraints on it.

\bigskip

Throughout this section, assume given a fixed symmetric-simplicial
group $G$ and work with a fixed term $G_n$ of $G$. 

\bigskip

The language and notations of section \ref{DoldKanDecomps} will be reused with one change as follows. A
{\it symmetric multi-index $\alpha$ of length $k$ and dimension $\le n$} 
is a strictly increasing sequence of indices
	$$\alpha = \big\{i_1 < \ldots < i_k \big\}$$
satisfying $1 \leq i_p \leq n$ for all $p$ (this is the change---the allowed range of the indices is shifted
upwards by 1). The length $\vert\alpha\vert$ 
of $\alpha$ is once again the number $k$ of indices. 
Denote the set of all symmetric multi-indices of dimension $\le n$ as follows.
	$$J^{(n)} := \big\{~\text{symmetric multi-indices $\alpha$ of dimension $\le n$}~\big\}$$
For brevity, symmetric multi-indices will be called simply multi-indices.  Hopefully this
will cause no confusion.

The ({\it generalized}) {\it quasidegeneracy operator} $u_\alpha$
corresponding to $\alpha$ is the composition of elementary quasidegeneracy
operators
	$$u_\alpha : G_{n-\vert\alpha\vert} \longrightarrow G_n$$
	$$u_\alpha := u_{i_k}u_{i_{k-1}}\ldots u_{i_1}$$
and there is a corresponding subgroup of quasidegenerate simplices
	$$U_\alpha := u_\alpha (M_{n-\vert\alpha\vert}) \subseteq G_n$$
for $\vert\alpha\vert\ge 1$. Note this subgroup is an embedded copy
of $M_{n-\vert\alpha\vert}$, since the quasidegeneracy operator $u_\alpha$ is
injective on account of the identity $d_iu_i = \id$ from Theorem \ref{QuasiIds}.
For $\alpha = \varnothing$, there is once again the corresponding subgroup
	$$U_\varnothing := M_n \lhd G_n.$$

\bigskip

As in section \ref{DoldKanDecomps}, let us consider the question of
which total orderings of the $U_\alpha$ give rise to SDP decompositions
of $G_n$ for any symmetric-simplicial group $G$. Such total
orders and their corresponding decompositions will be called {\it symmetric Dold-Kan}.

Our answer to this question is given in a similar manner. Letting
$J^{(n)}$ be partially ordered by inclusion, any total order of
$J^{(n)}$ extending the inclusion partial order will be a symmetric Dold-Kan
total order. 

\bigskip

In the rest of this section, $\alpha(\cdot)$ will denote a fixed total
order on the multi-indices of dimension $\le n$, thought of as a
bijection as follows.
	$$\alpha(\cdot): \big\{ 0, 1, \ldots, 2^n-1 \big\} \longrightarrow J^{(n)}$$
	$$k \mapsto \alpha(k)$$
On $\alpha(\cdot)$ the sole requirement is made that it extend the inclusion
partial order, or equivalently
	$$\alpha(k) \subseteq \alpha(l) \implies k \leq l.$$
Once again, this forces $\alpha(0)$ to be the empty multi-index
$\varnothing$ and $\alpha(2^n-1)$ to be the full multi-index.

\bigskip


The following theorem is the goal of the section. Write $U_k :=
U_{\alpha(k)}$.

\begin{thm} \label{quasimaindecomp}
	There is an internal SDP decomposition as follows.
		$$G_n = U_0 \rtimes U_1 \rtimes \ldots \rtimes U_{2^n-1}$$
\end{thm}

Here are tools analogous to those used in section \ref{DoldKanDecomps}. For an
arbitrary multi-index $\alpha$, the {\it face operator} corresponding in
a dual manner to $u_\alpha$ is
	$$d_\alpha : G_n \longrightarrow G_{n-\vert\alpha\vert}$$
	$$d_\alpha := d_{i_1}d_{i_2}\ldots d_{i_k}$$
There is a corresponding {\it projection operator} defined as follows.
	$$\pi_\alpha: G_n \longrightarrow G_n$$
	$$\pi_\alpha(g) := u_\alpha d_{\alpha}(g)$$
Similarly as in the statement of Theorem \ref{quasimaindecomp}, we again use
the notation $\pi_k := \pi_{\alpha(k)}$. By repeated
application of the identities $d_iu_i = \id$ (section \ref{SymmMooreComp}), one
obtains
	$$d_\alpha u_\alpha = \id$$
	$$\pi_k^2 = \pi_k.$$

The decomposition algorithm used here is identical in appearance to the one used
in section \ref{DoldKanDecomps}. Namely, for $g\in G_n$, define a sequence of elements
$g_i \in G_n$, starting with $g_{2^n-1} := g$ and proceed recursively in
reverse order by the formula 
	\begin{equation} \label{quasirecurse} g_{k-1} := g_k - \pi_k(g_k)
	\end{equation}
and ending with $g_0$. Thus one has the following decomposition for any
$g\in G_n$, in which one must once again be careful to keep the summands in
increasing order (left to right). \begin{equation}
\label{quasifactorize}
		g = g_{2^n-1} = \sum_{i = 0}^{2^n-1} \pi_i(g_i)
	\end{equation}
\bigskip

The proof of the SDP decomposition of Theorem \ref{quasimaindecomp}
above requires the following ``quasi'' analogs of Proposition \ref{idslemma} and
Theorem \ref{ids}. The proofs of these analogs are nearly the same as
those of the originals, so we leave the verifications to the interested
reader with the advice simply to follow the original proofs closely,
accounting carefully for the differences in the statements of
Propositions \ref{idslemma} and \ref{quasiidslemma}. and inserting the inclusion partial
order in place of the length-product partial order wherever it occurs.

\bigskip

The following convenient notation will be used.
	$$u_j^- := u_{j-1}$$

\begin{prp} \label{quasiidslemma} 
	Let $\beta = \big\{ j_1 < \ldots < j_l \big\}$ 
	be a symmetric multi-index of dimension $\le n$ and let $i$ be an index satsifying $1 \le i \le n$.
	Consider the following two alternatives, which are mutually exclusive and
	cover all possibilities for $\beta$ and $i$.  (The nonsensical statements
	$j_{l+1} > i$ and $i > j_0$ are agreed to be true for any $i$.)
	\begin{enumerate}
		\item{For some (unique)	subscript $0 \le q \le l$, one has $j_{q+1} > i > j_q$.}
		\item{For some (unique) subscript $1 \le q \le l$, one has $i = j_q$.}
	\end{enumerate}	
	In the first case, the identity 
			$$d_iu_\beta = u_{j_l}^-u_{j_{l-1}}^-\ldots 
				u_{j_{q+1}}^-u_{j_q}^{\phantom{-}}\ldots u_{j_1}^{\phantom{-}} d_{i-q}$$
	holds, and the rightmost factor $d_{i-q}$ is different from $d_0$.
	Then we say that $d_i$ has ``slipped past'' $u_\beta$. 
	In the second case, the identity 
			$$d_iu_\beta = u_{j_l}^-u_{j_{l-1}}^-\ldots 
				u_{j_{q+1}}^-u_{j_{q-1}}^{\phantom{-}}\ldots u_{j_1}^{\phantom{-}}$$
	holds, and we say that $d_i$ is ``absorbed by'' $u_\beta$. 
\end{prp}

\begin{thm} \label{quasiids}
	The following facts hold for arbitrary multi-indices 
		$$\alpha = \big\{ i_1 < \ldots < i_k \big\} ~~~~~~~~
		\beta = \big\{ j_1 < \ldots < j_l \big\}$$
		of dimension $\le n$.
	\begin{enumerate}
		\item{$d_{\alpha}u_\beta = u_{\beta^\prime}$ for some multi-index $\beta^\prime$
			if and only if $\alpha \subseteq \beta$.}
		\item{If $\alpha \nsubseteq \beta$ then $d_{\alpha}u_\beta 
				= d_{\alpha^\prime}u_{\beta^\prime}d_i$
				for some $\alpha^\prime,\beta^\prime$ and some $i \neq 0$.}
		\item{For $i \neq 0$, one has $d_id_{\beta} = d_{\beta^\prime}$
			for some multi-index $\beta^\prime \supsetneq \beta$.}
	\end{enumerate}
\end{thm}

\bigskip

The next three lemmas are the analogs of Lemmas \ref{kernel},
\ref{pikernel} and \ref{ids8}. The difference is that the the operators
$s_\alpha$ are replaced by $u_\alpha$ and the role of Theorem \ref{ids}
is taken over by Theorem \ref{quasiids}. The proofs are otherwise
completely identical, and so we omit them.

\begin{lem} \label{quasikernel}
	The following holds for any $k$ with $0 \leq k \leq 2^n-1$.
		$$g_k \in 
			\bigcap_{i = k+1}^{2^n-1}\mathrm{ker}~d_{\alpha(i)} = 
			\bigcap_{i = k+1}^{2^n-1}\mathrm{ker}~\pi_{i}$$
\end{lem}
\begin{lem} \label{quasipikernel}
	The following holds for any $k$ with $0 \leq k \leq 2^n-1$.
		$$\pi_k(g_k) \in U_k$$
\end{lem}
\begin{lem}\label{quasiids8}
	For any $k^\prime > k$, the operator $\pi_{k^\prime}$ annihilates $U_k$.
\end{lem}

\bigskip

The proof of Theorem \ref{quasimaindecomp} is now practically identical
to the proof of Theorem \ref{maindecomp}, and so we omit it as well.

\bigskip

One can define the length-product partial order for $J^{(n)}$
exactly as was done for $I^{(n)}$ in section \ref{DoldKanDecomps}.  The following
proposition shows that the length-product partial order extends the inclusion
partial order.

\begin{prp}\label{relatepos}
	For multi-indices $\alpha$ and $\beta$ as in Theorem \ref{quasiids}, 
	the following implication holds.
		$$\alpha \subseteq \beta \implies \alpha \prce \beta$$
\end{prp}
\prf{Proof.}{
	First note that 
		$$\alpha \subseteq \beta \implies \abs{\alpha} \le \abs{\beta}$$
	so that the 0th condition of $\alpha \prce \beta$ is satisfied.
	For the rest, define the following function.
		$$q(\cdot):\{1,\ldots,k\} \lra \{1,\ldots,l\}$$
		$$p \mapsto \text{~the~number~} q(p) \text{~such~that~} i_p = j_{q(p)}$$
	Note $q(\cdot)$ 
	is then a {\it strictly} increasing function, so that for each $p$ one has
		$$q(k-p) < q(k-p+1) < \ldots < q(k) \le l$$
		$$\implies q(k-p) \le l-p.$$
	Since the $j_q$ are increasing in $q$ one has immediately
		$$i_{k-p} = j_{q(k-p)} \le j_{l-p}$$
	as required in the definition of $\alpha \prce \beta$.
}

As a result of this proposition, any total order extending the
length-product partial order also extends the inclusion partial order.
In particular, the binary total order of section \ref{BinaryTO} can be
used as a symmetric Dold-Kan total order. This is done in the final
remark below in order to investigate the symmetric analog of the formula
from Remark \ref{CCremarkableform}.

\rmk{
	The following remark gives the symmetric analog of the formula of Carrasco and
	Cegarra from Remark \ref{CCremarkableform}.
	Fix an index $m$ with $0 \le m \le 2^n-1$, and let 
		$$\alpha = \{i_1 < \ldots < i_k\}$$ 
	stand for $\alpha(m)$. Also write $\alpha^c$ for the multi-index which
	is the complement of $\alpha$ in $\{1,\ldots,n\}$ and write
		$$\alpha^c = \{j_1 < \ldots < j_l\}.$$  
	Recall from Definition \ref{replacementop} the homomorphisms
		$$r_i: G_n \lra G_n$$
		$$r_i(g) := u_i d_i(g)$$
	for $1\le i \le n$, and let
		$$r_\alpha := r_{i_k} \ldots r_{i_1}$$
	be the product of $r_i$ for $i\in\alpha$ (order does not matter as they commute). 
	Also write $r_j^\perp$ for the crossed homomorphisms
		$$r_j^\perp: G_n \lra G_n$$
		$$r_j^\perp := 1 - r_j$$
		$$r_j^\perp(g) = g - u_jd_j(g)$$
	and write $r_{\alpha^c}^\perp$ for the product 
		$$r_{j_1}^\perp\ldots r_{j_l}^\perp$$
	of the $r_j^\perp$ for $j\in\alpha^c$ in {\it increasing} order
	from left to right.
\newline\indent	
	Then the component $\pi_k(g_k)$ of $g\in G_n$ under the SDP
	decomposition of the present section using binary order can be
	expressed explicitly in terms of $g$ by the following formula.
		$$\pi_k(g_k) = r_\alpha^{\phantom{\perp}} r_{\alpha^c}^\perp(g)$$
	Due to the relations
		$$r_i(1-r_j) = (1-r_j)r_i$$
	holding for all $i,j$, the factors $r_i$ may be mixed around among the factors $1-r_j$.
	In particular, they may be arranged into the product in increasing order
		$$\pi_k(g_k) = r_1^\eps r_2^\eps \ldots r_n^\eps$$
	where $r_i^\eps$ is $r_i$ or $(1-r_i)$ depending as $i$ does or does not belong to $\alpha$.
}